\documentclass[10pt]{amsart}
\author{Aram L. Karakhanyan}
\address{
School of  Mathematics\\
The University of Edinburgh\\
Mayfield Road, King's Buildings,  EH9 3JZ\\
Edinburgh, UK}
\email{aram.karakhanyan@ed.ac.uk}
\title[Homogenization]{On a conjecture of De Giorgi related to homogenization}

\author{Henrik Shahgholian}
\address{Department of Mathematics, KTH, Lindstedtsv\"agen 25, 100 44 Stockholm, Sweden}
\email{henriksh@kth.se}
\thanks{H. Shahgholian was supported by Swedish Research Council. A.Karakhanyan was partly supported by EPSRC grant. 
We thank Michael Benedicks for his insightful  comments on the dynamical system issues of the current note,
and Bj\"orn Engquist for bringing to our attention the paper \cite{Hou}.}

\usepackage{color}
\usepackage{graphics}
\usepackage[dvips]{graphicx}
\usepackage{amsmath}
\usepackage{amsfonts}
\usepackage{mathrsfs}
\usepackage{mathrsfs}
\usepackage{mathabx}                                                                                      %
\usepackage{upgreek}
\usepackage{esint} %%% \fint command for mean integral
\usepackage{fancyhdr}
\usepackage{latexsym}
\usepackage{amssymb}
\usepackage[all]{xy}
\usepackage{float,graphicx}
\usepackage{fancybox}
\usepackage{rotating}
\usepackage{tikz}
\usepackage{pgfplots}

\newtheorem{theorem}{Theorem}                                                                                %
\newtheorem{lemma}{Lemma}                                                                                    %
\newtheorem{proposition}{Proposition}                                                                        %
\newtheorem{definition}{Definition}                                                                          %

\theoremstyle{definition}                                                                                  %
\newtheorem{remark}[theorem]{Remark}

\theoremstyle{theorem}
\renewcommand\epsilon\varepsilon % use curly epsilon
\newcommand{\p}{\partial}

\newcommand\e{\varepsilon}

\newcommand\R{\mathbb R}
\newcommand\T{\mathbb T}
\newcommand\Z{\mathbb Z}

\newcommand\ol{\overline}

\renewcommand\phi{\varphi}

\newcommand\cal{\mathcal}

\newcommand{\lb}{\left(}
\newcommand{\rb}{\right)}

\newcommand\der[1]{\displaystyle\frac{d{#1}}{dt} }
\renewcommand{\nu}{\upnu}

\newcommand{\h}[1]{\textbf{#1}} % \bf-s the letters in math mode
\renewcommand\div{\operatorname{div}} % divergence
\newcommand\ds{\displaystyle}
\usepackage{amsmath}

\begin{document}

\begin{abstract} 
For a  periodic vector field $\h F$,  let
 $\h X^\e$ solve the dynamical  system 
\begin{equation*}
 \frac{d\h X^\e}{dt}  = \h F\lb\frac {\h X^\e}\e\rb .
\end{equation*}
In \cite{DeGiorgi}  Ennio De Giorgi enquiers  whether from the existence of the limit  
$\h X^0(t):=\lim\limits_{\e\to 0}\h X^\e(t)$ one can conclude that  $ \frac{d\h X^0}{dt}= constant$. 
Our main result settles this conjecture under fairly general assumptions on $\h F$, which in some cases may also depend on $t$-variable. 

Once the above problem is solved, one can apply the result to the corresponding transport equation, in a standard way. This is also touched upon in the text to follow.
\end{abstract}	

\maketitle

%\setcounter{tocdepth}{2}
%\tableofcontents

\section{Introduction}
\subsection{Problem setting}
For each $i=1, \dots, d$ let  $F_i:  [0,\infty) \times \R^d\to \R$ be  a smooth  1-periodic function in both variables. 
Let us  consider 
the first order  system of differential equations with oscillating structure 
\begin{equation}\label{eq-0eps}
  \der{x_i}=F_i\lb \frac{t}\e, \frac{x_1}\e, \dots, \frac{x_d}\e\rb\quad i=1, \dots, d, 
\end{equation}
where $\e>0$ is a small parameter. 
Our primary motivation for studying (\ref{eq-0eps}) comes from a conjecture 
posed by Ennio De Giorgi in
\cite{DeGiorgi} (Conjecture 1.1 page 175)
concerning the homogenization of the transport equation 
\begin{equation}\label{eq-0trans}
 \begin{array}{lll}
 \p_tu^\e( t, x )+ \h F\lb t/ \e, x/\e\rb  \cdot \nabla_x u^\e( t,x)=0, \quad t\in (0, \infty), x\in \R^d, \\
 u^\e(t=0, x)=u_0(x),  \quad \ x\in \R^d,
 \end{array}
\end{equation}
with vector field $\h F=(F_1, F_2, \dots, F_d)$ Lipschitz continuous and periodic in both variables $( t,x)$. The Lipschitz continuous initial condition  $u_0(x)$ is  specified at the initial time $t=0$.

\medskip 
He also conjectured that if (\ref{eq-0trans}) is homogenizable then the following property must be true (see 
\cite[ page 177]{DeGiorgi}):
{\it Let $\h X^\e(t)$ be the solution of the following initial value problem
\begin{equation}\label{ibvp}
 \der{\h X^\e}=\h F   \lb \frac{t}{\e},  \frac{ \h X^\e}\e\rb, \quad \h X^\e(0)= p
\end{equation}
for some given initial condition $p\in \R^d$.
Then the limit  exists
\begin{equation}\label{eff-lim}
 \h X^0(t):= \lim_{\e\to 0}\h X^\e(t)
\end{equation}
for any $t,  p$.
Moreover, there is a vector $\h B\in \R^d$ such that
\begin{equation}\label{eff-der}
 \frac{d \h X ^0}{d t}=\h B .
\end{equation}}
We remark that  Peirone \cite{Peirone}  showed that if $\h F$ does not depend on $t$ 
then the asymptotic linearity of $\h X^\e(t)$ as $t\to\infty$ implies that (2) is homogenizable, see Remark  \ref{semig}.

\subsection{Related work}
In view of Peirone's result \cite{Peirone}, the homogenization of (\ref{ibvp})
is closely related to the homogenization of the first order 
transport equations $\p_t u+\h F\cdot \nabla u=0$ describing miscible flow in porous media \cite{Tassa}. 
One of the central questions  concerning  (\ref{eq-0trans}) is the strong convergence which is not true 
in general as the example of  equation (\ref{eq-0trans}) with $\h F(t,x_1, x_2)=(0, \sin x_1), d=2$ shows,  see \cite{DeGiorgi} page 176.
It is known  that if
$\div \h F=0$
\footnote{This refers to the case of unit density $\rho=1$ for the invariant measure, see Section \ref{inv} for more details.} 
then the effective 
equation has arithmetic averages $(\int_{\T^2}F_1(x)dx, \int_{\T^2}F_2(x)dx)$ as the forcing velocity, 
whereas the shear field $\h F(x)=\h a \phi(x), \h a=(1,\gamma)\in \R^2$ yields harmonic averages, i.e. in the homogenized equation the forcing velocity is $\h a\int_{\T^2}\frac{dx}{\phi(x)}$, see \cite{Tassa}.
The interested reader can find more on this problem in the works  \cite{Tassa}, \cite{Regis} and \cite{Dali} and the references therein.

\medskip

The homogenization of more general transport equations  
\begin{equation}\label{gen-trnsp}
\p_t u^\e+\div[\h a_\e f(u^\e)]=0, \quad u^\e(0, x)=U_0(x, x/\e)
\end{equation}
under the assumption $\h a _\e =\h a (x, x/\e)$ and $\div_x\h a (x, y)=\div_y\h a(x, y)=0$, is studied in  \cite{Weinan}. 
The case when $\h a_\e=\h a(x/\e)$ is studied in \cite{Hou}. It is also shown that solutions of \eqref{gen-trnsp}
converge in $L^2$ and the limit equation is either a constant coefficient linear transport equation (ergodic case) 
or an infinite dimensional dynamical system, see  \cite{Weinan, Hou}.

In \cite{Tartar} Tartar studied some transport equations with memory effects. He addressed the question
of importance of considering the limit function rather than the equation it satisfies. The question he raised was
  whether the limit  retains, in some sense, the structure of linear transport equations  (e.g.,  when it is traveling wave solution).

Some of these questions were addressed by Tassa in \cite{Tassa}. In particular,  he showed that for shear flow  $(d=2)$ the limit is 
a traveling wave (Theorems 4.2 and 4.5 in \cite{Tassa}). He also  derived  convergence rate which depends on   the smoothness of the forcing vector field  as well as on 
whether the  rotation number (which we denoted $\gamma$ in the formula $\h a=(1, \gamma)$ above) is   rational or irrational. In fact for rational rotation number (Theorem 4.5 in \cite{Tassa}) 
the limit is determined by some function $a_\eta$ see (3.13) in \cite{Tassa}, and the limit function is a traveling wave if $a_\eta=const$ for all $\eta\in [0,1]$.

It seems plausible that the techniques here can (partially)
be applied to more general context involving random structure, i.e. stochastic differential equations.
Similar type of problems, have been studied in recent 
works of Bardi-Cesaroni-Scotti  \cite{BCS}. 
The problem here can be reduced to the well-known classical perturbation problem through variable substitution  $\h Y ^\e= \h X^\e /\e $. To illustrate this at a heuristic level, we assume (for clarity)  $\h F$ to be independent of $t$. We thus have 
$$ \der{\h X^\e}=\h F (\h Y^\e), \quad \h X^\e (0) = p\quad  \Rightarrow \quad  
 \der{\h Y^\e}=\frac{1}{\e} \h F (\h Y^\e), \quad \h Y^\e (0) = p/\e.$$ 
Introducing $\h Z^\e(s)=\h Y^\e(\e s), s>0$ we infer 
$\frac{d \h Z^\e}{ds}=\e  (\h Y^\e)'(\e s)=\h F(\h Z^\e)$ and thus 
by Theorem 3.1 \cite{Peirone} we get that for a fixed $\e>0$ the limit 
$$\lim_{s\to \infty}\frac{\h Z^\e(s)}{s}=\lim_{s\to \infty}\frac{\h X^\e(\e s)}{\e s}:=\h B$$
exists and is independent of $\e$ for a suitable class of $\h F$. 
If we knew that this limit is also \textit{uniform} in $\e$ then for $\tau=\e s$ we could conclude 
 that 
$\lim_{\e\to 0}\frac{\h X^\e(\tau)}{\tau}$ exists for each fixed $\tau$ and is independent of $\e$ or, equivalently, 
$\h X^\e(\tau)=o_\e(\tau)+\h B \tau$. Certainly this captures the case when $p=0$. Nevertheless, it is possible that our Theorem 2b has some overlapping with above mentioned Theorem 3.1 \cite{Peirone}. 
 
 A further direction, that our approach might be possible to extend to, is that of multi-scale problems. More exactly,
 one may consider ${\bf F}$ that has both slow and fast variable ${\bf F} (x,x/\e)$. A particular case of this was studied by 
 G. Menon \cite{Menon}, with ${\bf F} (x,x/\e) = \hbox{div} ({\bf K}(x) + \e {\bf A}(x/\e) ) $.

\subsection{Problem set-up} \label{subsec-rotnumb}
We shall switch between cases of $t$-dependent as well as $t$-independent $\h F$, and this will be clear from the context.  Hence we shall use both notation $\h F (t,x)$, as well as $\h F(x)$.

Next, going back to our $t$-independent $\h F$, one can establish  a number of remarkable properties, for the  non-oscillating system (i.e. when $\e=1$)
\begin{eqnarray}\label{eq-00}
  \der{x_i}=F_i\lb x_1, \dots, x_d \rb .
\end{eqnarray}
Suppose that (\ref{eq-00}) has invariant measure $d\mu_x=\rho(x)dx$ with density  $\rho >0$ i.e. the vector field $\rho  \h F$ is divergence free; see Section \ref{inv} for details. 
For the two dimensional problem,  $(d=2)$,  Kolmogorov 
proved that if $\h F(x)=(F_1(x_1, x_2), F_2(x_1, x_2))\not =0$, is $\Z^2$ periodic and  both 
$\rho $ and $\h F$ are real analytic in $(x_1,x_2)$ 
variables, then there is an analytic transformation of coordinates 
$y=\h f(x)$ such that (\ref{eq-00}) transforms into 
{\bf shear flow} system 
\begin{equation}\label{shear-eq-00}
\der{y_i}=\frac{a_i}{G(y_1, y_2)} , \quad i=1,2 ,
\end{equation} 
with constants $a_1=1, a_2=\gamma\in \R$ and $G$  being a  $\Z^2$ periodic scalar function. Here
$\gamma$ is called the {\bf rotation number of \eqref{eq-00}} (also called rotation index) and the system (\ref{eq-00}) is ergodic if $\gamma$ is diophantine, see \cite{Sinai}. For the latter case the 
shear flow \eqref{shear-eq-00} can be further transformed to a constant speed system 
$\der w^i=c_i, i=1, 2$ where $c_i$ are constants. 

In fact, one can take 
\[G(x_1, x_2)=\rho(x_1, x_2), \quad \mbox{for a.e.} \ \ (x_1, x_2)\in \R^2 \]
to be the density of invariant measure of \eqref{shear-eq-00} such that we have 
$\div \frac\rho G=0$. In other words, now $\frac1 G$ is the 
density of the invariant measure of the new shear flow system of differential equations \eqref{shear-eq-00}, obtained from \eqref{eq-00} via a coordinate transformation introduced by Kolmogorov \cite{Kolm53}.

\medskip

The main goal of this article is to analyze the behaviour of the solution 
$\h X^\e(t)$ to equation \eqref{ibvp}  as $\e\to 0$ under some conditions  imposed on the vector field
$\h F=(F_1, \dots, F_d )$ which we list below:

\begin{itemize}
\item[{\bf(F.1)}] $\h F:\R^d\to \R^d $ is continuous, $\Z^d$-periodic and there is a constant $L>0$ such that 
\begin{equation}
|\h F(\h u_1)-  \h F(\h u_2)|\leq L|\h u_1-\h u_2|, \quad \forall \h u_1, \h u_2 \in \R^d.
\end{equation}
We write $\h F=(F_1, \dots, F_d)$ where $F_i, 1\le i\le d$ are the components of the vector field $\h F$. 
\item[{\bf(F.2)}] There is a constant $\lambda>0$ such that 
 $$\lambda\le F_i(\h u)\leq \frac1\lambda, \quad 1\le i\le d$$ for every $\h u\in \R^d $. 
\item[{\bf(F.3)}]  There is a  bounded $\Z^d$ periodic  function $\rho>0$ such that $\div (\rho \h F)=0$ in $\R^d$. Here $\rho$ is called the density of invariant measure.
\end{itemize}

The equation $\div (\rho \h F)=0$ is understood in the  weak sense, i.e. 
$\int\rho \h F \cdot \nabla \psi=0$ for every $\psi\in C_0^\infty(\R^d)$.

The conditions $\bf  (F.1)-(F.3)$ will  be mainly used in the statement of Theorem 2.
%

%----------
% subsection
%----------

\subsection{Invariant measure}\label{inv} (General discussion)
Condition $\bf (F.3)$ needs some explanation. Suppose that 
$\h F(x)=\frac{\h a}{G(x)}, x\in \R^d$ for some constant vector $\h a$ and suitable scalar function $G$ such that $\h F$ is smooth. It is clear that for this case $\rho=G$. However for general flows the existence of $\rho$ is not easily obtained. In the proof of Theorem 2b below 
we require that the invariant measure exists and is  bounded in order to construct a change of variables which reduces general flows to shear one.
In this regard we mention the following existence result from \cite{Lions}:
Suppose $\h F :\R^d\to \R^d$, $\h F\in C^1$  and for simplicity $t$-independent.
Let $\rho$ be sought as the solution of Liouville's equation
$\div(\rho\h F)=0$. Let $\tau=x_d, x'=(x_1, x_2, \dots, x_{d-1}, 0)$ and assume that $F_d>0$ then Liouville's equation can be rewritten as follows
$$\partial_\tau\log \rho+\nabla_{x'}\log\rho\cdot \frac{\h F'}{F_d}=-\frac{\div \h F}{F_d}  ,  $$
where $\h F'=(F_1, \dots, F_{d-1}, 0)$.
We can specify initial condition at time $\tau=0$, i.e. $x_d=0$ and then by  \cite{Lions} (Proposition II.1 and Remark afterwards)  
there is a $L^\infty$-solution 
of this Cauchy problem in $\R^d\times[0, \infty)$, provided that both $\h F$ and  the initial data are Lipschitz.

If $d=2$ then it is well known that divergence free vector field is 90 degree rotation of  the gradient of a potential function $u$, i.e.,  $\rho \h F=(u_{x_2}, -u_{x_1})$. From here  we have that 
$\rho F_1=u_{x_2}, \rho F_2=-u_{x_1}$. For $\h F$ satisfying $(\h F.2)$ we can eliminate $\rho$ to obtain 
\begin{equation*}
u_{x_1}=-u_{x_2}\frac{F_2}{F_1} \quad \mbox{in}\ \T^2.
\end{equation*}
The existence and regularity of periodic solution 
 $u=u(x_1, x_2)$ follows from standard existence theory
for the first order linear equations via the method of characteristics. In particular if $\h F\in C^k$ then $\nabla u\in C^k$. 
The density of the invariant measure can be recovered as follows $\rho=\frac{u_{x_1}F_2-u_{x_2}F_1}{|\h F|^2}$.

%----------
% subsection
%----------

\subsection{The approach and methodology}

 De Giorgi's conjecture has (more or less)  been ignored completely. 
Indeed,   the fact that convergence of the underlying dynamical system would give the 
convergence of the transport problem,  have been unnoticed in the literature.  Our result (read observation) should be seen in the light of  homogenization of the dynamical system,
 rather than the transport problem; even though this directly implies the convergence of the transport problem. 
The approach we have taken here is a combination of a few, already worked out, methods (originating in the work of Kolmogorov \cite{Kolm53}, and later Bogolyubov \cite{Bogolyubov}). More precisely, it is a combination of Kolmogorov's transformation of coordinate system (and its refinement due to Tassa \cite{Tassa}) and Bogolyubov's method for singular perturbations. In particular, the latter implies a convergence rate 
as $\e\to 0$.

To the best of our knowledge, this has not been done previously and hence worth noticing.
Such a   composition of hybrid techniques -- of combining singular perturbations, dynamical systems and homogenization --  
gives new insights and opens up for the study of convergence rates for similar problems. 

We also want to stress that although our result seems to be new,   it does not use any new technique, and most probably if the problem was noticed by others, that have worked with the related transport problem,  a similar observation would have been made.

\medskip

%----------
% section
%----------
\section{Preliminaries and main results}

We first recall the definition of KBM-functions from \cite{Sanders} Definition 4.2.4.

\begin{definition}\label{def-KBM}
Consider the function $G(t, x)$ 
continuous in $t$ and $x$ on $[0,\infty)\times \R^d$
such that for some constant $L>0$
there holds
\[|G(t, x_1)-G(t, x_2)|\le L|x_1 - x_2|, \quad \hbox{for all}\ t\in[0, \infty), x_1, x_2\in \R^d.\]
If the average 
\begin{equation}\label{eq-def-G0}
G^0(y)=\lim_{\ell\to \infty}\frac1 \ell\int_0^\ell  G(\tau , y)d\tau 
\end{equation}
exists uniformly in $y$ on compact sets $D\subset \R^d$
then we call $G$ a KBM-function (KMB stands for Krylov, Bogolyubov and Mitropolski.)
\end{definition}

\medskip 
We next justify the existence of $G^0$ and obtain a refined estimate 
for $\delta$ under the periodicity assumption on $G$ in $t$-variable.

\begin{lemma}\label{lem-blya}
Consider the function 
$G: [0, \infty)\times \R^d\to \R, $ continuous in $t\in [0,\infty)$ and $x\in \R$ 
such that for some constant $L>0$
there holds
\[|G(t, x_1)-G(t, x_2)|\le L|x_1 - x_2|, \quad \hbox{for all}\ t\in[0, \infty), x_1, x_2\in \R^d.\] 
Suppose  $G(t, x)$ is  $1$-periodic in $t$,  then the limit in \eqref{eq-def-G0} exists and consequently $G$ is a KBM-function. 

\end{lemma}

\begin{proof}
For fixed $y$ we have  
\begin{eqnarray*}  
\int_0^\ell(G(\tau , y)-G^0(y))d \tau&=&\sum_{m=1}^{[\ell ]}\int_{m-1}^{m}(G(\tau, y)-G^0(y))d\tau+\int_{[\ell]}^\ell (G(\tau, y)-G^0(y))d\tau=\\\nonumber
&=&\int_{[\ell]}^\ell (G(\tau, y)-G^0(y))d\tau
\end{eqnarray*}
where $[\ell]$ is the integer part of $\ell>0$ and $G^0(y)=\int_0^1G(\tau, y)d\tau$.
Consequently, 
\[\frac1{\ell}\int_0^\ell G(\tau, y)d\tau=G^0(y)+\frac1\ell\int_{[\ell]}^\ell (G(\tau, y)-G^0(y))d\tau\to G^0(y)\  \hbox{as}\ \ell\to \infty.\]
The second part that $G$ is KBM follows from Lemma 4.6.4 \cite{Sanders}. 

\end{proof}

Note that for periodic $G$ independent of $x$, we have that 
$G^0(x)=\int_0^1 G(\tau,x)d\tau$ is constant.
The convergence rate for almost periodic $G$ depending only on $x$ variable may be weaker  as the example in 
Section \ref{sec-exmpl} shows.
\medskip

We formulate our main results below starting from the one dimensional problem.

\begin{theorem} ($d=1$)
 Let  $G_\star(t,x), t\in \R, x\in \T$ be positive, periodic in $x$,  such that the function $G(t, x):=G_\star(x, t)$ (with swapped variables) is KBM-function and  
 $$M:=\sup_{x\in \R, t\ge0}\frac1{G_\star(t, x)}<\infty.$$

 Let $X^\e$ be the solution to the
initial value problem
\begin{eqnarray}\label{eq-Th1}
 \left\{
\begin{array}{ccc}
\ds \der{X}^\e=\frac1{ G_\star\lb t, \frac{ X^\e}\e \rb },\\
X^\e(0)=p.
\end{array}
\right.
\end{eqnarray}
Then there is a Lipschitz continuous  function $X^0(t)$ such that 
\begin{equation}\label{Th1-a}
| X^\e(t)-X^0(t)|\le C(T)\e, \quad t\in[0, T]
\end{equation}
where $T>0$ is the length of the time interval $t\in [0, T]$, $C(T)$ is a positive constant depending only on $T$ and $G_\star$.
Furthermore, if  $G_\star(t, \eta)$ does not depend on $t$ and 
is  periodic in $\eta$, then $X^0(t)=p_0+\beta t$ for some $p_0, \beta\in \R.$ 
\end{theorem}

In the proof of Theorem 1 we will use a simple version of Bogolyubov's method, 
tailored for the Cauchy problem $\der{\h Y^\e} =\h G\lb\frac t\e , \h Y^\e\rb$, $\h Y^\e(0)=p$, 
see   \cite{Bogolyubov} \S 26, \cite{Sanders} Lemma 4.3.1. 
It is worthwhile to mention that at some point we swap the arguments 
of the function $G_\star$ such that the resulted function $G$ is KBM.

\medskip 

Next we state our main result for the multidimensional problem.

\begin{theorem}\label{thm-d-dim} ($d\geq 2$)
\smallskip 
\begin{itemize}
\item[\bf a)] Let $d\ge 2$ 
be a periodic scalar function $G: \R^d\to \R$ independent of $t$, 
$G\in C^k(\R^n)$ and there are positive constants  $c_0$ and $\kappa, k>d+\kappa+1$ 
such that $\h a\in \R^d$ is diophantine, i.e., 
 \begin{eqnarray}\label{eq-dphnt}
|\langle \h a, \h m\rangle |\geq \frac {c_0}{|\h m|^{d+\kappa}},\quad \forall \h m\in \Z^d\setminus\{0\}.
\end{eqnarray}
Finally, suppose that 
$\h F=\frac{\h a}G$ satisfies {\bf (F.1)-(F.3)}. 
   If $\h X^\e$ is the solution to 
the Cauchy problem $\der{\h z^\e}=\frac{\h a }{G\lb \frac{\h z^\e}{\e}\rb}, \h z^\e(0)=p$  then
$$
\left|{\h z^\e(t)} -\left(p+\frac{\h a }{\cal M(G)}t\right)\right| \leq C \e,\quad  t\ge 0
$$
where  $\cal M(G)=\fint_{\T^d} G$, $\T^d$ is the $d$-dimensional torus,  and 
\[C=%\frac{d|\h a |}{\pi |\cal M(G)|}\sum\limits_{m\in\Z^d\setminus\{0\}}\frac{|G_m|}{ |\langle m, \h a\rangle| }\le 
\frac{ d|\h a|}{c_0\pi |\cal M(G)|}\sum\limits_{m\in\Z^d\setminus\{0\}}| m|^{d+\kappa}{|G_m|}{ }<\infty\]
 where $G_m$ is the $m$-th Fourier coefficient of $G$.

\item[\bf b)] Let $d=2$ and $\h F  \in C^k$ be   independent of $t$,  and 1-periodic in $x$-variable.
Let further {\bf (F.1)-(F.3)} hold and $\h X^\e$ solves the Cauchy problem 
\begin{eqnarray}\label{eq-Th2}
 \left\{
\begin{array}{ccc}
\ds \der{\h X}^\e=\h F\lb \frac{\h X^\e}\e \rb, \\
\h X^\e(0)=p.
\end{array}
\right.
\end{eqnarray}
 Let $\gamma$ be the rotation number (see section \ref{subsec-rotnumb}) and assume that $\h a =(1, \gamma)$ satisfies \eqref{eq-dphnt} with some constants $C>0$ and $\kappa>0$ such that 
$k>3+\kappa$. 
Then there is a linear function $\h X^0(t)=p+\h B t$, $\h B\in \R^2$ such that
\begin{equation}\label{Th1-b}
|\h X^\e(t)-\h X^0(t)|\le \hat C\e, \quad t\in[0, \infty)
\end{equation}
where $\hat C$ depends on $ \|\rho\h F\|_\infty$, $\gamma$ and $\|\h F\|_{C^k}$.
\end{itemize}
\end{theorem}

We shall use a number of results from dynamical systems.  In particular, in the proof of Theorem 2 we shall employ Kolmogorov's theorem on coordinate transformation  $y=\h f(x)$ \cite{Kolm53}, see section \ref{subsec-rotnumb}. It needs to be mentioned that Kolmogorov's proof is not constructive i.e., he did not write 
explicit form of such transformation. In \cite{Tassa} Tassa found a simple argument that renders the explicit form of $\h f$.
Such coordinate transformation exists for $d\ge 3$ under various assumptions   \cite{Arn92}, \cite{Ko07}.

%\medskip 
%%%%%%%%%%%%%%%%%%%%%%%%%%%%%%%%%%%                                                                                                           
%                                                                                                          %
%                                           Section                                                   %
%                                                                                                          %
%%%%%%%%%%%%%%%%%%%%%%%%%%%%%%%%%%%
\section{Proof of Theorem 1}

We first observe that if $G$ is a KBM-function then by Definition \ref{def-KBM} the following limit 
\begin{equation}\label{G0}
G^0(y)=\lim_{\ell\to \infty}\frac1 \ell\int_0^\ell  G(s, y)ds
\end{equation}
exists uniformly in $y\in D$   for any compact $D\subset \R$.
In particular, the Lipschitz continuity of $G$ translates to $G^0$.
Next let us derive a    scaled  version of  Bogolyubov's estimate in one dimension.

\begin{lemma}\label{lem-tech}
Let $G: [0, \infty)\times\R\to \R$ be a KBM-function periodic in the first variable.
Let $h^\e(\xi)$ be the solution of the Cauchy problem 
$\displaystyle \frac{d h^\e}{d\xi}=G\lb \frac \xi \e, h^\e(\xi)\rb$, $h^\e(0)=p$.

Let $G^0$ be as in \eqref{G0} and $h^0$ a unique solution of the Cauchy problem 
$$ \frac{d h^0}{d\xi}= G^0(h^0), \qquad  h^0(0)=p$$
on the finite interval $[0, T_1]$. Then, as $\e\to 0$, 
\begin{equation}\label{rate-1}
| h^\e (\xi)-h^0  (\xi)|\le C(T)\e , \qquad  0\leq  \xi \leq T <T_1
\end{equation}
for some constant $C(T)>0$ depending only on $T$.
\end{lemma}

\begin{remark}
Note that under the conditions of Theorem 1 the solution $h^0$ is unique because $G^0$ is Lipschitz.
\end{remark}

\begin{proof}
We use  Bogolyubov's estimate for the slowly varying systems.
Define $\theta^\e(\xi)={h^\e(\e \xi)}$ 
then we have 
\begin{eqnarray}
\left\{
\begin{array}{ccc}
\ds \frac{{d\theta}^\e}{d\xi}= \e G\lb \xi, {\theta^\e} \rb, \\
\theta^\e(0)=p.
\end{array}
\right.
\end{eqnarray}
Furthermore, let $\theta^0(\xi)$ solve 
\begin{eqnarray}
\left\{
\begin{array}{ccc}
\ds \frac{{d \theta}^0}{d\xi}= \e G^0\lb {\theta^0} \rb, \\
\theta^0(0)=p ,
\end{array}
\right.
\end{eqnarray}
where $G^0$ is as in \eqref{G0}. 
Applying Bogolyubov's estimate,  \cite{Chechkin} Theorem 12.1 and Remark 12.1, (see also \cite{Sanders} Theorem 4.5.5)  to 
$\theta^\e, \theta^0$ we have that 
\begin{equation}\label{eq-Bog}
\sup _{\xi\in[0, \frac T\e]}|\theta^\e(\xi)-\theta^0(\xi)|\leq C(T)\e
\end{equation}
where 
$C(T)>0$ depends only on $T$.
After setting $h^0(\e \xi)=\theta^0(\xi)$, substituting $\e \xi =s$ in \eqref{eq-Bog}  
the result follows.
 \end{proof}

\medskip

Now we are ready to finish the proof of Theorem 1.
Observe that  $0< \der{X^\e}\le \sup_{x\in \R, t\ge 0} \frac1{G_\star(t, x)}=M<\infty$ and therefore $\{X^\e\}$ is uniformly Lipschitz continuous on every 
finite interval $[0, T]$. In fact, we  have the estimate $|X^\e(t)|\le |p|+TM, t\in [0, T]$. Furthermore, $X^\e$ is strictly monotone because  $G>0$. 
Thus $X^\e$ has inverse which we denote by 
$h^\e$,
 \begin{equation}\label{X-inv}
\xi = X^\epsilon (h^\e(\xi)).
\end{equation}

Rewriting the system for $h^\e$ we have
$$
\frac1{\frac{d h^\e}{d\xi }}=\frac1{G_\star(h^\e(\xi), \xi/\e)}\qquad \Rightarrow \qquad {\frac{d h^\e}{d\xi }}={G_\star(h^\e(\xi),\xi/\e)}.
$$
As for the initial condition, we have   $h^\e(p)=0$.

Denote $G(t, x)={G_\star(x, t)}$, the function with swapped variables. Note 
that $G$ satisfies all requirements of  Lemma \ref{lem-tech} (in particular $G$ is periodic in $t$), and  hence it follows that  $h^\e\to h^0$ locally uniformly on 
$[0, \infty)$ and the homogenized equation is $\frac{dh^0}{d\xi}=G^0(h^0)$
where $$G^0(y)=\lim_{\ell\to \infty}\frac1 \ell\int_0^\ell{G(\tau,y)}d\tau.$$

Returning to $X^\e$ and using the refined convergence rate (\ref{rate-1}) for  periodic $G$ in $t$ variable, we note that by \eqref{X-inv}
\begin{equation}
\xi =X^\e\left(h^\e(\xi)\right)=X^\e([h^\e(\xi)-h^0(\xi)]+h^0(\xi))  
\end{equation}
implying that $|\xi-X^\e(h^0(\xi))|=|X^\e([h^\e(\xi)-h^0(\xi)]+h^0(\xi))  -X^\e(h^0(\xi))  |\le M|h^\e(\xi)-h^0(\xi)|\le MC(T)\e$, where the last inequality follows from \eqref{rate-1}. Hence, $X^\e$ converges uniformly to $X^0(t)$, determined by the implicit equation $\xi=X^0(h^0(\xi))$.

Finally, the last part of Theorem 1 follows from the fact that  $G^0$
is constant for  periodic $G_\star$ and therefore $X^0 (t) $ must be linear function of $t$.

%%%%%%%%%%%%%%%%%%%%%%%%%%%%%%%%%%%                                                                                                           
%                                                                                                          %
%                                           Section                                                   %
%                                                                                                          %
%%%%%%%%%%%%%%%%%%%%%%%%%%%%%%%%%%%

\medskip 

\section{Multi-dimensional problem: Proof of Theorem $2a$}

\subsection{Change of variables for $d\ge 2$} 
Let $d\mu=\rho dx $ be the invariant measure of the system $\der{ x_i}=F_i(x)$ 
where $\h F(x)=(F_1(x), \dots, F_d(x)), x\in \R^d$ is the vector field on the right hand side of the  equation
(\ref{eq-00}).	If $d=2$, $F_1, F_2, \rho\in C^\infty$, $F_i:\R^2\to \R$ and $F_1^2+F_2^2>0$ then Kolmogorov showed that 
there is a transformation  $x\to y$ such that in the new system of coordinates 
the equation transforms into the shear flow $\der {y_1}=F, \der {y_2}=\gamma F$	
where $\gamma$ is the rotation number (see section \eqref{subsec-rotnumb}), and $F$ is a positive function.
Furthermore, if $\gamma$ is diophantine (see the formulation of Theorem 2 for precise condition) then there is another transformation of $\R^2$, $y\to u$ such that
the system takes the form $\der {u_i}=a_i, i=1,2$ where $a_i$ are constants.

\smallskip 

For $d\ge 3$ Kolmogorov's theorem has been generalized by   Kozlov which we state below without proof, see \cite{Ko07}. 

\begin{proposition}\label{prop-1}
Let $d\ge 2$ and $G>0, \frac1 G\in C^{k}, k>d+\kappa+1$ is smooth. If $\h a=(a_1, \dots, a_d)$ is diophantine in the sense of \eqref{eq-dphnt}
then there exists a change of variables transforming
the system 

\begin{eqnarray}\label{shear-00}
 \der{w_j}=\frac{a_j}{G(w_1, \dots, w_d)}, \qquad j=1, \dots, d
\end{eqnarray}
 into the constant coefficient system
$\der{w_j}=a_j$.
\end{proposition}

\smallskip

It is clear that for the shear flow \eqref{shear-00} the density of invariant measure is 
$\rho=G$.

\medskip 

\subsection{Proof of Theorem 2a}

\begin{proof}
We shall use the coordinate transformation introduced in \cite{Ko07} Theorem 2: if  $\h u(t)$ solves the 
shear system $\der{\h u}=\frac{\h q}{G(\h u)}$ with diophantine $\h q$ then the mapping given by 
the equations
\begin{equation}\label{blya-1}
w_j=u_j+\frac{q_j}{{\mathcal M(G)}}f(\h u), \qquad 1\le j\le d, \qquad  \h u =(u_1, \dots, u_d)
\end{equation}
transforms the equation into 
$$\der{w_j}=\frac{q_j}{\mathcal M(G)  },$$
as stated in Proposition \ref{prop-1}, see \cite{Ko07} page 197. Here $\mathcal M(G)=\fint\limits_{\T^d}G$ is the mean value of $G$ and $f$ is determined from the 
first order differential equation 
$$\langle \nabla f, \h q\rangle=G(\h u)-\mathcal M(G ).$$ 
In fact, this mapping is non-degenerate (i.e. has nontrivial Jacobian) 
and  is  one-to-one \cite{Ko07}. 
Taking $\e \h u= \h z$ we see that 
$$\der{\h u}=\frac{\h q}{G(\h u)}$$ 
with $\h q=\frac{\h a}\e$. 
From Fourier's expansion we have
$$G(\h u)-\mathcal M(G)=\sum\limits_{m\in\Z^d\setminus\{0\}}G_me^{2\pi i\langle m, \h u\rangle}$$ 
which by integration  gives 
\begin{equation}\label{sill}
f(\h u)=\sum\limits_{m\in\Z^d\setminus\{0\}}\frac{G_m}{2\pi i\langle m, \h q\rangle }e^{2\pi i\langle m,\h  u\rangle} =\e \sum\limits_{m\in\Z^d\setminus\{0\}}\frac{G_m}{2\pi i\langle m, \h a\rangle }e^{2\pi i\langle m,\h  u\rangle} ,
\end{equation}
and  the series is absolutely convergent,  due  to the assumption that $\h a $ is diophantine (see \eqref{eq-dphnt}) and $G\in C^k, k>d+\kappa+1$. In particular,  $|G_m|\le C(k, d)(1+|m|)^{-k}$ for some universal constant $C(k, d)>0$ depending only on $d$ and $k$. Notice that 
the sum is bounded because $\frac 1G$ satisfies the assumptions $\bf (F.1)-(F.3)$. 
Summarizing we have 
 
\begin{eqnarray}\label{blya-2}
w^\e_j(t)
&=&\frac{q_j}{\cal M(G)} t+w^\e_j(0) \\\nonumber&=&
\frac{a_j}{\e \cal M(G)} t+w^\e_j(0).
\end{eqnarray}
On the other hand from \eqref{blya-1} and \eqref{sill}
\begin{eqnarray}\label{blya-3}
w^\e_j(t)
&=&
u^\e_j(t)+\frac{q_j}{\cal M(G)}f(\h u^\e)\\\nonumber
&=&
\frac{z_j^\e(t)}\e+ \left\{\frac{a_j}{\e} \frac1{\mathcal M(G)}\right\} \e \sum\limits_{m\in\Z^d\setminus\{0\}}\frac{G_m}{2\pi i\langle m, \h a\rangle }e^{2\pi i\langle m, \h u^\e\rangle}\\\nonumber
&=&
\frac{z_j^\e(t)}\e+ \frac{a_j}{\mathcal M(G)}\sum\limits_{m\in\Z^d\setminus\{0\}}\frac{G_m}{2\pi i\langle m, \h a\rangle }e^{\frac{2\pi i}\e \langle m, \h  z^\e\rangle}.
\end{eqnarray}
Combining \eqref{blya-2}, \eqref{blya-3} and $w_j^\e(0)=\frac{z_j^\e(0)}\e+\frac{q_j}{\mathcal M(G)}f(\frac{z^\e(0)}\e)$, which follows from \eqref{blya-1}, 
we get 
$\frac{a_j}{\cal M(G)}t+p_j-z_j^\e(t)=\cal \sigma(\e)$
where 
\begin{eqnarray*}
\mathcal \sigma(\e)&=&\e\left\{-\frac{q_j}{\mathcal M(G)}f(\frac{p}\e)+ \frac{a_j}{\mathcal M(G)}\sum\limits_{m\in\Z^d\setminus\{0\}}\frac{G_m}{2\pi i\langle m, \h a\rangle }e^{\frac{2\pi i}\e \langle m, \h  z^\e\rangle}\right\}\\
&=&
\e\left\{\frac{a_j}{\mathcal M(G)}\sum\limits_{m\in\Z^d\setminus\{0\}}\frac{G_m}{2\pi i\langle m, \h a\rangle }\left[e^{\frac{2\pi i}\e \langle m, \h  z^\e\rangle}-e^{\frac{2\pi i}\e \langle m,  p\rangle}\right]\right\}.
\end{eqnarray*}
Since $G\in C^k, k>d+\kappa+1$ and $\h a$ is diophantine, see \eqref{eq-dphnt}, it follows that 
the series 
$\sum\limits_{m\in\Z^d\setminus\{0\}}\frac{|G_m|}{2\pi |\langle m, \h a\rangle| }$ converges. 
Therefore using \eqref{eq-dphnt} \[|\mathcal \sigma (\e)|\le \frac{2\e |\h a|}{|\cal M(G)|}\sum\limits_{m\in\Z^d\setminus\{0\}}\frac{|G_m|}{2\pi |\langle m, \h a\rangle| }\le 
\frac{\e |\h a|}{c_0\pi |\cal M(G)|}\sum\limits_{m\in\Z^d\setminus\{0\}}| m|^{d+\kappa}{|G_m|}{ }
\] and the series converges because from $G\in C^k$ we get $|G_m|\le C(k, d)(1+|m|)^{-k}$ with $k>d+\kappa+1$. The proof now follows.
\end{proof}

\medskip

\begin{remark}\label{semig}
Peirone showed that if $\h F\in C^1(\T^d)$ is $\Z^d$ periodic, $u_0\in C^1$ 
 and the limit $\lim\limits_{t\to \infty} \frac{S^t_{\h F}(x)}{t}$ exists 
for a.e. $x\in\T^d$  then the problem 
(\ref{eq-0trans}) is homogenizable, see \cite{Peirone} Lemma 2.2 (b).
Here $S_{\h F}^t$ is the semigroup generated by (\ref{eq-00}).
Our result establishes the converse of this statement for homogenizable (\ref{eq-0trans}).
\end{remark}

%%%%%%%%%%%%%%%%%%%%%%%%%%%%%%%%%%%                                                                                                           
%                                                                                                          %
%                                           Section                                                   %
%                                                                                                          %
%%%%%%%%%%%%%%%%%%%%%%%%%%%%%%%%%%%

\section{Proof of Theorem $2b$}

Our goal here is to apply Kolmogorov's coordinate transformation in order to reduce the general problem to shear flow.
For this, Tassa \cite{Tassa} found an explicit formula,  that we will write below. 
We should (again)  point out that Kolmogorov's proof in \cite{Kolm53} is not constructive.   

It is convenient to introduce some basic facts about the equation $\der{\bf X}=\h F(\h X)$ with $\h F$ satisfying the 
properties ({\bf F.1})-({\bf F.3}). Let $d\mu=\rho dx$ be the invariant measure corresponding to this system, 
then by definition $\div(\rho \h F)=0$. Thus the vector field $\h b=(b_1,b_2)=\rho\h F$
is divergence free, 1-periodic, and $\rho\in C^k,$ see section \ref{inv}. This yields that  the integral $\int_0^1 b_1(x_1,x_2)dx_2$ is constant since
\begin{eqnarray}
\p_{x_1}\int_0^1 b_1(x_1,x_2)dx_2&=&\int_0^1\p_{x_1} b_1(x_1,x_2)dx_2\\\nonumber
&=&
-\int_0^1 \p_{x_2}b_2(x_1,x_2)dx_2\\\nonumber
&=&-[b_2(x_1,1)-b_2(x_1,0)]=0.
\end{eqnarray}
Similarly we have that $\int_0^1 b_2(x_1,x_2)dx_1$ is constant. Denote
$\ol{b_1}=\int_0^1 b_1(x_1,x_2)dx_2, $ $\ol{b_2}=\int_0^1 b_2(x_1,x_2)dx_1$ (which are 
the mean integrals of $b_1, b_2$ over $\T^2$) and set
\begin{eqnarray}\label{y-coord}
\begin{array}{lll}
\displaystyle y_1=f_1(x_1,x_2)=\frac1{\ol{b_2}}\int_0^{x_1}b_2(\xi,0)d\xi,\\
\vspace{-0.3cm}\\
\displaystyle y_2=f_2(x_1,x_2)=\frac1{\ol{b_1}}\int_0^{x_2}b_1(x_1,\xi)d\xi.
\end{array}
\end{eqnarray}
It is shown in \cite{Tassa} that in the new coordinate system we get the 
shear flow
 $\der{\h y}=\frac{\h a}{G(\h y)}$ with $\h a=(1,\gamma)$,  where $\gamma$ is   the
rotation number, see  section \ref{subsec-rotnumb}. Furthermore, we have that 
\begin{equation}\label{determinant}
\left|\frac{\p( y_1, y_2)}{\p( x_1, x_2)}\right|=\frac{b_1(x_1, x_2)}{\ol{b_1}} \frac{b_2(x_1,0)}{\ol{b_2}}\not=0, \quad \forall x\in \T^2
\end{equation}
and the invariant measure density is 
\begin{equation}
\frac1{G(y)}=\frac{b_2(g_1(y),0)}{\ol {b_2}}F_1(g_1(y), g_2(y))
\end{equation}
with $\h g=(g_1, g_2)$ being the inverse of $\h f=(f_1, f_2)$, see \cite{Tassa}, page 1395.   
In particular, it follows \begin{equation}\label{ref-bbb}
\h b\in C^k\quad \hbox{and}\quad G\in C^k
\end{equation} 
(recall that $\h b=\rho \h F$ and $\rho\in C^k,$ see section \ref{inv}). Moreover from  \eqref{y-coord}, \eqref{determinant} and the inverse function theorem $\h g\in C^k$ implying $G\in C^k, k>3+\kappa$.
\medskip

In order to take advantage of (\ref{y-coord}) we introduce the function $\h z^\e(t)=\h X^\e(t)/\e$. Then $\h z^\e(t)$ solves the Cauchy problem 
$\der{\bf z^\e}=\frac{\h F(\h z^\e)}\e$,  $\h z^\e({0})=\frac p\e$. Clearly, the invariant measure now is  
$d\mu_z=\frac1\e\rho dz$ and $\h b^\e=(\frac{b_1}\e,\frac{b_2}\e)$ is divergence free. Note that 
$$\frac{b_i}{\ol{b_i}}=\frac{b_i^\e}{\ol{b_i^\e}}$$ and therefore applying the change of variables 
$y=\h f(x)$, with mapping $\h f=(f_1, f_2)$ given by (\ref{y-coord}) we obtain the shear flow
\begin{equation}
\der{\h y^\e}=\h a \frac{b_2(g_1(\h y^\e), 0)}{\ol{b_2}}\frac{F_1(\h g(\h y^\e))}{\e}.
\end{equation}

In order to get rid of $\e$ in the denominator  we set $\h w^\e(t)=\e \h y^\e(t)$. Then $\h w^\e(t)$ solves the equation 
\begin{equation}
\der{\h w^\e}=\h a \frac{b_2(g_1(\h w^\e/\e), 0)}{\ol{b_2}} {F_1(\h g(\h w^\e/\e))}.
\end{equation} 

By \eqref{y-coord} we have that 
\begin{eqnarray*}
f_1(x_1+1, x_2)&=&f_1(x_1, x_2)+1,\\\nonumber
 f_1(x_1, x_2+1)&=&f_1(x_1, x_2),
\end{eqnarray*}
and similarly 
\begin{eqnarray*}
f_2(x_1+1, x_2)&=&f_2(x_1, x_2),\\\nonumber
 f_2(x_1, x_2+1)&=&f_2(x_1, x_2)+1
\end{eqnarray*}
in view of the periodicity of $\h b$.  
Consequently if $\h e_i$, $(i= 1,  2)$ is the unit vector in the canonical basis of $\R^2$ then this translates to the 
inverse of $\h f$, namely we have  $g_j(\eta+\h e_i)=g_j(\eta)+M_{ij}, 1\le i,j\le 2$
where $M_{ij}\in \Z$, see \cite{Tassa}  equation  (2.5). 
This yields that $\frac1{G(\eta)}=\frac{b_2(g_1(\eta), 0)}{\ol{b_2}}{F_1(\h g(\eta))}$ is periodic 
function and $\h w^\e$ solves the Cauchy problem 
\begin{equation}
\der{\h w^\e}=\frac{\h a} {G\lb\frac{\h w^\e}\e\rb}, \quad  \h w^\e(0)=\e \h f\lb \frac{\h x^\e(0)}\e\rb.
\end{equation}
From here, in light of  (\ref{y-coord}) we have 
\begin{eqnarray}
w_1^\e(t)&=&\frac\e{\ol{b_2}}\int_0^{z_1^\e(t)}b_2(\xi, 0)d\xi=\frac\e{\ol{b_2}}\int_0^{x_1^\e(t)/\e}b_2(\xi, 0)d\xi\\\nonumber
&=& x_1^\e+\frac\e{\ol{b_2}}\int_{0}^{x_1^\e(t)/\e} (b_2(\xi, 0)-\ol{b_2})d\xi\\\nonumber
&=& x_1^\e+\frac\e{\ol{b_2}}\int_{[x_1^\e(t)/\e]}^{x_1^\e(t)/\e} (b_2(\xi, 0)-\ol{b_2})d\xi\\\nonumber
\end{eqnarray}
for  $\h b$ is periodic, see  the proof of Lamma \ref{lem-blya} for a similar argument. Here $[\cdot]$ denotes the integer part.

 Hence we conclude that 
\begin{equation}\label{x1-asym}
w_1^\e(t)=x_1^\e(t)+
\frac\e{\ol{b_2}}\int_{[x_1^\e(t)/\e]}^{x_1^\e(t)/\e} (b_2(\xi, 0)-\ol{b_2})d\xi
\quad t\in [0, \infty).
\end{equation}
In particular for the initial condition we get 
that $w_1^\e(0)=p_1+\cal O(\e)$.
As for the asymptotic expansion of $w_2^\e$ then we need to use a well-known fact that 
there is a scalar function $\phi$ such that $\h b=(\p_2\phi, -\p_1\phi)$ for every 
two dimensional divergence free vector field $\h b\in L^\infty$. From this equation it follows 
$\phi(x)=\psi(x)+\h q\cdot x+q_0$ where $\psi$ is periodic. 
Observe that $\h b=\rho \h F\in L^\infty$ by $\bf (F.1)-(F.3)$ (in fact $b\in C^k$ by \eqref{ref-bbb}),  hence 
$\|\nabla \psi\|_\infty\le \|\rho\h F\|_\infty<\infty$. Using this fact we compute
\begin{eqnarray*}
w_2^\e(t)&=&\frac\e{\ol{b_1}}\int_0^{z_2^\e(t)}b_1\lb \frac{x_1^\e}\e, \xi \rb d\xi =\frac\e{\ol{b_1}}\int_0^{x_2^\e(t)/\e}b_1\lb \frac{x_1^\e}\e, \xi\rb d\xi\\
&=& \frac\e{\ol{b_1}}\int_0^{x_2^\e(t)/\e}\left[b_1\lb \frac{x_1^\e}\e, \xi\rb -b_1(0, \xi)\right]d\xi+\frac\e{\ol{b_1}}\int_0^{x_2^\e(t)/\e}b_1\lb 0, \xi\rb d\xi\\
&=&  
\frac\e{\ol{b_1}} \left[\phi \lb \frac{x_1^\e}\e, 0\rb-\phi \lb \frac{x_1^\e}\e, \frac{x_2^\e}\e\rb +\phi \lb 0, \frac{x_2^\e}\e\rb +\phi(0, 0)\right]
+\frac\e{\ol{b_1}}\int_0^{x_2^\e(t)/\e}b_1\lb 0, \xi\rb d\xi\\
&=&  
\frac\e{\ol{b_1}} \left[\psi \lb \frac{x_1^\e}\e, 0\rb-\psi \lb \frac{x_1^\e}\e, \frac{x_2^\e}\e\rb +\psi \lb 0, \frac{x_2^\e}\e\rb +\psi(0, 0)+2q_0\right]\\\nonumber
&&+\  x_2^\e(t)+\frac\e{\ol{b_1}} \int_{[x_2^\e(t)/\e]}^{x_2^\e(t)/\e} (b_1\lb 0, \xi\rb -\ol{b_1})d\xi\\
\end{eqnarray*}
where the third line follows as in (\ref{x1-asym}), or integrating by parts and using $\h b=(\p_2\phi, -\p_1\phi)$.
In particular, at $t=0$ we have that $w_2^\e(0)=
p_2+\cal O(\e )$.

Summarizing, we see that $\h w^\e$ solves the following Cauchy problem
\begin{equation*}
\der{\h w^\e}=\frac{\h a}{G\lb\frac{\h w^\e}\e\rb}, \qquad  \h w^\e(0)=\lb p_1+\cal O(\e),p_2+\cal O(\e)\rb
\end{equation*}
where $\h a=(1, \gamma)$ and $\gamma$ is the rotation number, see section \ref{subsec-rotnumb}.
By Theorem  2a, for diophantine $\h a=(1, \gamma)$ 
there is a linear function $\h w^0$ such that 
$|\h w^\e(t)-\h w^0(t)|\le C \e, t\in[0, \infty)$. Here $C>0$ depends on $\h F, \rho$ and $\gamma$ as in Theorem 2a (note that we can apply Theorem 2a because by \eqref{ref-bbb} $G\in C^k$). 
Then from \eqref{x1-asym}
$$|x_1^\e(t)-w_1^0(t)|\leq \frac\e{\ol{b_2}}\left|\int_{[x_1^\e(t)/\e]}^{x_1^\e(t)/\e} (b_2(\xi, 0)-\ol{b_2})d\xi\right| \le \frac{2\e\|\rho \h F\|_{L^\infty}}{\ol{b_2}}. $$
Finally for $x_2^\e$ we have
\begin{eqnarray*}
|x_2^\e-w_2^0|&\le &|w_2^\e - w_2^0|+\e\|\psi\|_{L^\infty}+\frac\e{\ol{b_1}}\left| \int_{[x_2^\e(t)/\e]}^{x_2^\e(t)/\e} (b_1\lb 0, \xi\rb -\ol{b_1})d\xi\right|\\
  &=&\hat C\e+\e\|\psi\|_{L^\infty}+\e \frac{2\|\rho\h F\|_{L^\infty}}{\ol{b_1}}
  \end{eqnarray*}
and the desired estimate follows

%%%%%%%%%%%%%%%%%%%%%%%%%%%%%%%%%%%                                                                                                           
%                                                                                                          %
%                                           Section                                                   %
%                                                                                                          %
%%%%%%%%%%%%%%%%%%%%%%%%%%%%%%%%%%%
\medskip

\section{Examples}\label{sec-exmpl}
\noindent
{\bf Example 1:}
Let $\h F$ be 1-periodic vector field such that $F_2=1$ and 
$$ F_1(x_1, x_2)=F_1(x_1)=\left\{ 
\begin{array}{ll}
1 &\qquad  0<x_1 \leq 1/2 ,\cr
0 &\qquad   1/2 < x_1 \leq 1.
\end{array}
\right.
$$ 
Let $\h X^\e(t)$ be the solution to the following initial value problem
$$
 \left\{ 
\begin{array}{ll}
\der{\h X^\e} =\h F\lb \frac{\h X^\e}\e \rb,\cr
\h X^\e(0)=p. 
\end{array}
\right.
$$
Let $S_\e$ be the  (cigar-shaped)  $\frac{\e}{2\sqrt{2}}$-neighborhood of the ray $p+sE, E=(2,1), s\ge 0$, i.e.  
$$S_\e=\left\{x\in \R^2: |x-[p+sE]|\le \frac{\e}{2\sqrt{2}}, s\ge 0\right\}.$$
Thus as $\e\to 0$ the trajectory (i.e. curves determined by $\h X^\e$) converges to the line $\ell(s)=p+s(2,1), s\ge 0$ in Hausdorff distance.
Hence the trajectory of the limit is the line $\ell(s)$.
As for the speed of the convergence, we note first that by definition $x^\e_2=1$ and it is enough to study the ode
$\der z=F_1(z/\e)$. Multiplying both sides of this equation by $\der{z^\e}$ and integrating 
we obtain that 
\begin{eqnarray*}
\int_{0}^s\left|\der{z^\e(t)}\right|^2dt=\int_0^s F_1 \lb \frac{z^\e(t)}\e  \rb \der{z^\e(t)}dt=qs+\cal O(\e)
\end{eqnarray*}
where $q=\fint_{[0, 1]}F_1=\frac12$. Since $\left|\der{z^\e(t)}\right|\leq \sup F_1=1$ we can use a customary compactness argument and infer from Lebesgue's dominated 
 convergence theorem 
 \begin{eqnarray*}
\int_{0}^s \left|\der{z^0(t)}\right|^2dt=\frac s2,
\end{eqnarray*}
where $z^0$ is the limit function. After differentiation we get $\left|\der{z^0(t)}\right|=\sqrt{ \frac12}$.

The astute reader has probably noticed that we did not use condition $\bf (F.3)$ here, but could still obtain a convergence rate. This is due to the one-dimensional character of the problem, since $F_2=1$ here.

\medskip

%%%%%%%%%%%%%%%%%%%%%%%%%%%%%%%%%%%

\noindent
{\bf Example 2:}   (One-dimension) Another example is given by $F$ with  saw-like graph
\begin{eqnarray}
F(\tau)=\left\{
\begin{array}{ccc}
\frac{2h\tau}a +\sigma& {\rm if} \ \tau\in [0, \frac a2),\\
\frac{2h}a(a-\tau)  +\sigma& {\rm if}\ \tau\in[\frac a2, a),
\end{array}
\right.
\end{eqnarray}
periodically extended over $\R$, see Figure 1. Here $\sigma>0$, $a>0$ is the periodicity of $F$
and $h=\max F$ is the peak of $F$. We can solve this equations explicitly: 
indeed we have that
\begin{eqnarray*}
\der{y^\e}=\left\{
\begin{array}{ccc}
\frac{2h}a(\frac{y^\e}\e-ka) +\sigma& {\rm if} \ \frac{y^\e}\e \in ak+[0, \frac a2),\\
\frac{2h}a(a(k+1)-\frac{y^\e}\e)  +\sigma& {\rm if}\ \frac{y^\e}\e \in ak+[\frac a2, a).
\end{array}
\right.
\end{eqnarray*}
After integration one gets
\begin{eqnarray*}
{y^\e}=\left\{
\begin{array}{lll}
C_-(k)e^{\frac{2ht}{\e a}}+\e ka -\frac{a\sigma\e}{2h}& {\rm if} \ {y^\e} \in \e ak+[0, \frac{ a\e}2),\\
C_+(k)e^{-\frac{2ht}{\e a}}+\e (k+1)a -\frac{a\sigma\e}{2h}& {\rm if} \ {y^\e} \in \e ak+[\frac{ a\e}2, a\e),
\end{array}
\right.
\end{eqnarray*}
with some constants $C_\pm(k)$ and $k\in \Z$. 
Clearly this solution $y^\e$ is monotone and hence the argument using the inverse function in the proof of Lemma 
\ref{lem-tech} works here too.
Obviously $\frac1a\int_0^a\frac{d\tau}{ F(\tau)}=\frac a h\log\lb \frac{ h+\sigma}{\sigma} \rb\equiv\beta$ and therefore we infer that 
$y^\e$ converges uniformly to $y^0(t)=p+\frac t{\beta}$ on any finite closed interval $[0, T]$.

 \begin{center}
 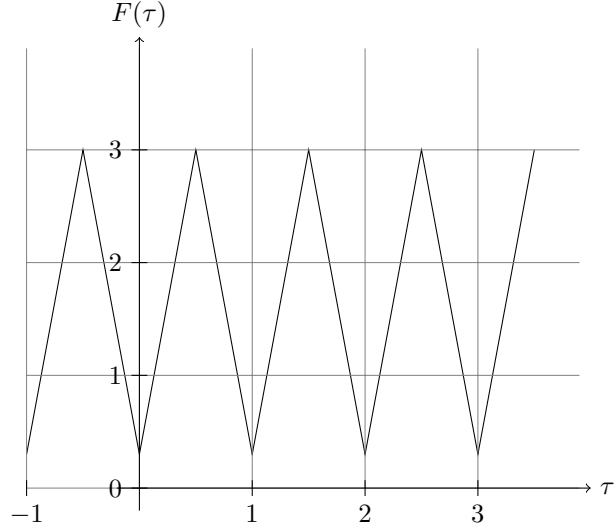
\begin{figure}
 \begin{tikzpicture}[scale=1.5]
\draw (-1,0.3) -- (-0.5,3) --(0,0.3) --(0.5,3) -- (1,0.3) -- (1.5,3) -- (2,0.3) -- (2.5,3) -- (3,0.3) --(3.5,3);

  \draw[style=help lines] (-1,0) grid (3.9,3.9)
       [step=0.25cm]      (1,2) 
              ;

  \draw[->] (-0.2,0) -- (4,0) node[right] {$\tau$};
  \draw[->] (0,-0.2) -- (0,4) node[above] {$F(\tau)$};

  \foreach \x/\xtext in {-1/-1,1/1,  
  2/2, 3/3}
    \draw[shift={(\x,0)}] (0pt,2pt) -- (0pt,-2pt) node[below] {$\xtext$};

  \foreach \y/\ytext in { 0/0, 1/1, 2/2, 
  3/3}    
    \draw[shift={(0,\y)}] (2pt,0pt) -- (-2pt,0pt) node[left] {$\ytext$};

 \end{tikzpicture}
 \caption{In this example $a=1, h=3$ and $\sigma=\frac13$.}
 \end{figure}
 \end{center}
\medskip

 \noindent
 {\bf Example 3:} Theorem 1 is still valid if the periodicity of $G(\cdot, x)$  is 
 replaced with almost periodicity in $x$ because  we  needed  periodicity in the proof only for  the convergence 
 rate for $G^0$. In this case one may get weaker error estimates, see \cite{Vu} Example 11.13. Indeed, 
 the function  $F(x)=\sum_{k=0}^\infty\frac1{(2k+1)^2}\sin\lb \frac x{2k+1}\rb$ is almost periodic. 
 By direct computation 
 
 \begin{eqnarray*}
 \int_0^TF(x)dx&=&\sum_{k=0}^\infty\frac2{2k+1}\sin^2\lb \frac T{2(2k+1)}\rb\\
 &=&\sum_0^{N(\e)}\dots+\sum_{N(\e)}^\infty\dots
 \end{eqnarray*}
 $N(\e)\sim \frac1\e$ then in this case $|X^\e(t)-X^0(t)|\le C(T)(\e|\log\e|)$
 on finite time intervals $[0, T]$.

\medskip 

\noindent
{\bf Example 4:}  (1-dimensional Transport Equation) One can apply  Theorem 1 to the homogenization of  some model transport equations 
such as 
\begin{equation}\label{eq-trans1D}
\p_t v^\e+H(x/\e)\p_x v^\e=0,\quad v^\e( 0,x)=v_0(x).
\end{equation}
Here $H>0$ is $C^1$ smooth periodic function.  Let $\h F=(1, H)$ and $\rho$ be the density of invariant measure, i.e.
$\div(\rho \h F)=0$. Therefore there is a function $M( t,x)$ solving the system 
\begin{eqnarray}
\left\{
\begin{array}{ccc}
\p_t M= -\rho H(y/\e)\\
\p_x M=\rho.
\end{array}
\right.
\end{eqnarray} 
The level sets $M=const$ are the characteristics of the equation \eqref{eq-trans1D}.
Noting that $\p_t M=-\rho H\not =0$ and applying the inverse function theorem to $M( t,x)=const$
we infer that $x=h^\e(t)$ and therefore for the solution of the Cauchy problem we have the 
formula 
$$v^\e( t,x)=v_0(x-h^\e(t)),$$
where by construction $\der{h^\e}=H(h^\e/\e), h^\e(0)=x$. Denote $v^0( t,x)=v_0(x-h^0(t))$, where $h^0= \lim h^\e$.  Thus we have from Theorem 1
the estimate 
\begin{eqnarray*}
|v^\e(t,x)-v^0(t,x)|&=& |v_0(x-h^\e(t))-v_0(x-h^0(t))|\\
&\leq&
 \|\p_x v_0\|_\infty|h^\e(t)-h^0(t)|\le C(T) \e
 \end{eqnarray*}
on finite time intervals $[0, T].$
\medskip 

\bibliographystyle{plain}

\begin{thebibliography}{1}

\bibitem{Arn92}
V.I. Arnold,
\newblock Polyintegrable flows.
\newblock {\em Algebra i Analiz}, 4(6):54--62, 1992.

\bibitem{BCS} M. Bardi, A. Cesaroni, A. Scotti,  Convergence in Multiscale Financial Models with Non-Gaussian Stochastic Volatility, preprint.


\bibitem{Bogolyubov} N. N. Bogolyubov, 
Y. A. Mitropolski, Asymptotic methods in the theory of non-linear oscillations, translated from Russian, New York: Gordon and Breach, 1961 

\bibitem{Chechkin}  G.A. Chechkin,  A.L. Piatnitski, A.S. Shamaev,  Homogenization: Methods and Applications (Translations of Mathematical Monographs), AMS 2007

\bibitem{Dali}
A. L. Dalibard, Homogenization of linear transport equations in a stationary ergodic setting, Comm. Partial Differential Equations, 33 (2008), pp. 881--921.


\bibitem{DeGiorgi} E. De Giorgi, On the convergence of solutions of some evolution
differential equations,
Set-Valued Analysis
1994, Volume 2, Issue 1-2, pp 175-182.

\bibitem{Lions} R.J. DiPerna, P.-L. Lions, 
Ordinary differential equations, transport theory and Sobolev spaces. 
Invent. Math. 98 (1989), no. 3, 511--547




\bibitem{Weinan} W. E, Homogenization of linear and nonlinear transport equations. 
Comm. Pure Appl. Math. 45 (1992), no. 3, 301--326




\bibitem{Hou} T. Hou, X. Xin, Homogenization of linear transport equations with oscillatory vector fields. 
SIAM J. Appl. Math. 52 (1992), no. 1, 34--45.


\bibitem{Regis} H. Ibrahim,  R. Monneau, 
On the Rate of Convergence in Periodic Homogenization of Scalar First-Order Ordinary Differential Equation, 
SIAM J. Math. Anal., 42(5), 2155--2176.




\bibitem{Kolm53}
A.N. Kolmogorov,
\newblock On dynamical systems with an integral invariant on the torus.
\newblock {\em Doklady}, 93(5):763--766, 1953.

\bibitem{Ko07}
V.V. Kozlov,
\newblock Dynamical systems with multivalued integrals on a torus.
\newblock {\em Proceedings of Steklov institute of Mathematics},
  256(1):188--205, 2007.



\bibitem{Menon} G. Menon,  Gradient systems with wiggly energies and related averaging problems. 
Arch. Ration. Mech. Anal. 162 (2002), no. 3, 193Ð246. 

\bibitem{Peirone} R. Peirone, Convergence of solutions of linear transport equations, 
Ergod. Th. and Dynam. Sys. (2003), 23, 919--933.

\bibitem{Picci} L. Piccinini, Homogeneization problems for ordinary differential equations, Rend. Circ. Mat. Palermo (2), 27 (1978), pp. 95--112.


\bibitem{Sanders} J. Sanders, F. Verhulst, J. Murdock, 
Averaging Methods in Nonlinear Dynamical Systems,
 Springer 2007.


\bibitem{Sinai} Ya. Sinai,
Introduction to ergodic theory,  Princeton University Press, Princeton, N.J., 1976.

\bibitem{Tartar} L. Tartar, Nonlocal effects induced by homogenization, in PDE and 
Calculus of Variations, pp 925--938, F.Culumbini et al., eds., Birkh\"auser, Boston, 1989.

\bibitem{Tassa} T. Tassa, 
Homogenization of two-dimensional linear flows with integral invariance. 
SIAM J. Appl. Math. 57 (1997), no. 5, 1390--1405. 

\bibitem{Vu} F. Verhulst,
Methods and Applications of Singular Perturbations, Texts in Applied Mathematics 50, Springer, 2010.
\end{thebibliography}

\end{document}